\documentclass{elsarticle}

\usepackage{lineno,hyperref}
\modulolinenumbers[5]

\journal{Journal of XXXXX}
\usepackage{amssymb,amsmath,amsfonts}
\usepackage{graphicx}
\usepackage{graphics}
\usepackage{epstopdf}
\usepackage[latin1]{inputenc}
\usepackage{url}
\newcommand{\thmlist}{
\begin{list}{Step 1}
{\setlength{\leftmargin}{0.6 in}\setlength{\labelwidth} {0.5 in}}}
\newcommand{\alglist}{
\begin{list}{Step 1}
{\setlength{\leftmargin}{1.1 in} \setlength{\labelwidth}{1.0 in}}}
 \newcommand{\proof} {\noindent {\bf Proof.} \quad}
 \newcommand{\eproof} {$\quad \square$}
 
 \newtheorem{theorem}{Theorem}[section]

 \newtheorem{definition}{Definition}[section]
 \newtheorem{property}{Property}[section]
 \newtheorem{remark}{Remark}[section]
\newcommand{\subtitle}[1]{\color{blue}}

\begin{document}

\begin{frontmatter}

\title{Symplectic Geometric Methods for Matrix Differential
Equations Arising from Inertial Navigation Problems}

\author{
Xin-Long Luo \textsuperscript{a}
and Geng Sun \textsuperscript{b} \\
\textsuperscript{a}School of Information and Communication Engineering, \\
Beijing University of Posts and Telecommunications, P. O. Box 101, \\
Xitucheng Road  No. 10, Haidian District, Beijing China 100876, \\
Tel.: +086-13641229062, luoxinlong@bupt.edu.cn, \\
https://teacher.bupt.edu.cn/luoxinlong; \\
\textsuperscript{b}Institute of Mathematics,
Academy of Mathematics and Systems Science, \\
Chinese Academy of Sciences,
sung@amss.ac.cn
}

\begin{abstract}
  This article explores some geometric and algebraic properties of
  the dynamical system which is represented by matrix differential
  equations arising from inertial navigation problems, such as the symplecticity
  and the orthogonality. Furthermore, it extends the applicable fields of
  symplectic geometric algorithms from the even dimensional Hamiltonian system to
  the odd dimensional dynamical system. Finally, some numerical experiments are
  presented and illustrate the theoretical results of this paper.
\end{abstract}


\end{frontmatter}


\section{Introduction}

Simultaneous localization and mapping (SLAM), which is involved in the inertial
navigation and robotic vision, is a hot research topic in some engineering
fields for which there are many applications such as unmanned aerial vehicles,
autonomous vehicles, remote medical operations and so on
(see \cite{Corke2016, LM2013, QLS2018, SMP2017}). In the front
end of SLAM, in order to track the pose of camera which is fixed to the
vehicle, it needs to solve a matrix differential equation arising from the
inertial navigation problem. In engineering fields, they usually adopt the
explicit Runge-Kutta method for this differential equation
owing to its simplicity of implementation (see \cite{Corke2016, QLS2018}).

\vskip 2mm

It is well-known that the explicit Runge-Kutta method is not
suitable for a Hamiltonian dynamical system, since the non-symplectic Runge-Kutta
method can not preserve its geometric and algebraic properties such
as the symplecticity and the orthogonality \cite{FQ2010, GSS2013,
HLW2006, Higham1997, LR2004, Sun1993, Sun2000, SGLS2018, WX2013, XT2004}.
Hong, Liu and Sun \cite{HLS2005} consider the symplecticity of a Hamiltonian
system which is represented by a PDEs with the skew-symmetric matrix
coefficient.

\vskip 2mm

For the inertial navigation problems, the variables are
a $3 \times 3$ matrix. We known that the dimension of variables is odd.
Thus, it can not directly apply the classical results of the symplectic geometric
algorithm to this problem, since the classical symplectic geometric algorithm
is applicable to the even dimensional Hamiltonian system. On the other hand,
it is important to preserve the geometric and algebraic properties of the
differential equation for the numerical method. Therefore, in this article,
we investigate some geometric and algebraic properties of the matrix differential
equation such as the law of generalized energy, the pseudo-symplecticity and
the orthogonal invariant. Consequently, we discuss the geometric and algebraic
properties of the symplectic Rune-Kutta method for the linear matrix
differential equation with the skew-symmetric matrix coefficient in Section 3.
In Section 4, we compare the symplectic Runge-Kutta method with the non-symplectic
Runge-Kutta for the differential equation with the skew-symmetric matrix
coefficient. The simulation results illustrate the theoretical results of this
article. Finally, in Section 5, we give some conclusions and discuss the
future work.

\vskip 2mm

\section{Geometric Structure of Linear Matrix Differential Equations}

\vskip 2mm

We choose a moving coordinate system connected to the aerial vehicle and
consider motions of the aerial vehicle where the origin is fixed. By
one of Euler's famous theorems, any such motion is infinitesimally a rotation
around an axis. We represent the rotational axis of the aerial vehicle by the
direction of a vector $\omega$ and the speed of rotation by the length of
$\omega$. Thus, the velocity of a mass point $x$ of the aerial vehicle is given
by the exterior product
\begin{align}
  v = \omega \times x = \begin{bmatrix}
  0 & - \omega_3 & \omega_2 \\
  \omega_3 & 0 & -\omega_1 \\
  -\omega_2 & \omega_1 & 0
  \end{bmatrix}
  \begin{bmatrix}
  x_1 \\
  x_2 \\
  x_3 \\
  \end{bmatrix},
\end{align} \label{RELSPEED}
which is orthogonal to $\omega$, orthogonal to $x$, and of length
$\|\omega\|\|x\| \sin\gamma, \; \cos\gamma = \frac{x^{T}\omega}{\|\omega\| \|x\|}$.

\vskip 2mm

We also regard the motion of the aerial vehicle from a coordinate system
stationary in the space. The transformation of a vector $x\in R^{3}$ in the
aerial vehicle frame, to the corresponding $y \in R^{3}$ in the stationary
frame, is denoted by
\begin{align}
  y = R(t) x. \label{MOVSYSTOSTASYS}
\end{align}
Matrix $R(t)$ is orthogonal and describes the rotation of the aerial vehicle. For
$x = e_i$ in the aerial vehicle frame, we find that the columns of matrix $R(t)$
are the coordinates of the aerial vehicle's principal axes in the stationary
frame. From equation \eqref{RELSPEED}, we know that these rotate with the
velocity
\begin{align}
  \left(\omega \times e_1,\; \omega \times e_2, \; \omega \times e_3 \right)
  = \begin{bmatrix}
  0 & -\omega_3 & \omega_2 \\
  \omega_3 & 0 & -\omega_1 \\
  -\omega_2 & \omega_1 & 0
  \end{bmatrix} = W. \label{PASPRE}
\end{align}
Thus, we obtain $\dot{R}$, which is the rotational velocity expressed in the
stationary frame, by the back transformation \eqref{MOVSYSTOSTASYS}:
\begin{align}
  R^{T} \dot{R} = W, \nonumber
\end{align}
which gives
\begin{align}
  \dot{R} = RW,  \hskip 2mm  R(t_0)^{T}R(t_0) = I,    \label{LMDER}
\end{align}
where $S$ is a skew-symmetric matrix, namely $W^{T} = -W$, and $R$ represents
the rotational transform matrix (see p. 48 in \cite{Corke2016}, p. 122
in \cite{MP2017}, or p. 278 in \cite{HLW2006}). If we denote $Q(t) = R(t)^{T}$
and $S = W^{T}$,
from equation \eqref{LMDER}, we have
\begin{align}
  \dot{Q} = SQ,  \hskip 2mm  Q(t_0)^{T}Q(t_0) = I,    \label{LMDE}
\end{align}
where $S$ is a skew-symmetric matrix. We denote the solution of
linear differential equation \eqref{LMDE} as
\begin{align}
 \Phi_{t}(Q(t_0)) = Q(t;\; Q(t_0)), \label{FLOWQ}
\end{align}
and term the map $\Phi_{t}: \; \mathbb{R}^{M \times M } \to \mathbb{R}^{M \times
M}$ as the flow map of the given system \eqref{LMDE}.

\vskip 2mm

We introduce some concepts before we investigate geometric structures of
equation \eqref{LMDE}. Assume that $f: \mathbb{R}^{M \times M} \to \mathbb{R}$
is a smooth function. Its directional derivative along a matrix
$X \in \mathbb{R}^{N \times N}$
is denoted here by
\begin{align}
  df(X) = \lim_{\Delta t \to 0} \frac{f(Z + \Delta t X) - f(Z)}{\Delta t}
  = \sum_{i = 1}^{N}\sum_{j = 1}^{N} \frac{\partial f}{\partial z_{ij}}x_{ij},
  \label{DIRDERI}
\end{align}
where the partial derivatives of $f$ are computed at a fixed location
$Z \in \mathbb{R}^{M \times M}$, and $z_{ij}, \; x_{ij}$ are the elements
of matrix $Z$ and matrix $X$, respectively. The linear function $df(\cdot)$ is called the
\textit{differential} of $f$ at $Z$ and is an example of a differential one-form.

\vskip 2mm

Using this denotation, we define the wedge product $df \wedge dg$ of two
differentials $df$ and $dg$  as follows (see pp. 61-64 in \cite{LR2004}):
\begin{align}
\left(df \wedge dg\right)(X, \, Y) = df(X)dg(Y) - df(Y)dg(X), \label{BLINFUN}
\end{align}
where $X, \; Y \in \mathbb{R}^{M \times M}$. Thus, we give the definition of
the wedge product of two matrix functions $P \in \mathbb{R}^{M \times M}$ and
$Q \in \mathbb{R}^{M \times M}$
as follows (see p. 30, \cite{Tao2006}):
\begin{align}
 dP \wedge dQ =
  = \sum_{i=1}^{M} \sum_{j = 1}^{M} \left(dp_{ij} \wedge dq_{ij}\right),
 \label{SMFORM}
\end{align}
where $p_{ij}$ and $q_{ij}$ are the elements of matrix $P$ and $Q$,
respectively. For the wedge product \eqref{SMFORM} of two matrix differentials
$dQ$ and $dP$, it also has some basic properties similar to the wedge product of
two vector differentials $dp$ and $dq$. We state them as the following Property
\ref{PROMWEDG}.

\begin{property} \label{PROMWEDG}
Let $dP, \, dQ, \, dR$ be $M \times M$-matrices of differential one-forms on
$\mathbb{R}^{M \times M}$, then the following properties hold.
\begin{itemize}
  \item[1.] Skew-symmetry
  \begin{align}
   dP \wedge dQ &= - dQ \wedge dP.  \label{SKESYM}
  \end{align}
  \item[2.] Bilinearity: for any $\alpha, \, \beta \in \mathbb{R}$,
  \begin{align}
   dP \wedge (\alpha dQ + \beta dR) &
 = \alpha dP \wedge dQ + \beta dP \wedge dR. \label{BILINE}
  \end{align}
  \item[3.] Rule of matrix multiplication
  \begin{align}
   dP \wedge (AdQ) & = (A^{T}dP) \wedge dQ, \label{MATMUL}
  \end{align}
  for any $M \times M$ matrix $A$.
  \end{itemize}
\end{property}

\vskip 2mm

Now we give the pseudo-symplectic property of equation \eqref{LMDE}.

\begin{property} \label{PSESYM}
Assume that $Q(t)$ is the solution of equation \eqref{LMDE}
with an initial orthogonal matrix $Q(t_0)$, then it satisfies the following
geometric property
\begin{align}
  dQ \wedge S dQ  = const,  \label{PSLDM}
\end{align}
where the wedge product $dP \wedge dQ$ is defined by equation \eqref{SMFORM}
and $const$ is a constant number.
\end{property}

\vskip 2mm

\proof Actually, if $Q(t)$ is the solution of equation (\ref{LMDE}), from properties
\eqref{SKESYM}-\eqref{MATMUL},  we have
\begin{align}
\frac{d}{dt}(dQ \wedge S dQ) & = d\dot{Q} \wedge S dQ + dQ \wedge S d\dot{Q}
 = SdQ \wedge SdQ + S^{T}dQ \wedge d\dot{Q} \nonumber \\
 & = -SdQ \wedge d\dot{Q}  = -SdQ \wedge SdQ  = 0, \nonumber
\end{align}
which proves the result of equation \eqref{PSLDM}. \eproof

\vskip 2mm
\begin{remark}If the skew-symmetric matrix $S$ equals
$J =\begin{bmatrix} 0 & I_{d}
\\ -I_{d} & 0 \end{bmatrix}$ and $Q \in \mathbb{R}^{2d}$ is the vector
function, equation \eqref{PSLDM} is the condition of symplecticity
(see p. 65 in \cite{LR2004}).
\end{remark}

\begin{theorem} (Kang Feng \& Zai-Jiu Shang \cite{FS1995}).
Assume that $Q(t)$ is the solution of the matrix differential equation
\eqref{LMDE} and matrix $S$ is skew-symmetric. Then, the determinant $det(Q(t))$
is an invariant.
\end{theorem}
\proof We denote $g(t) = det(Q(t))$, then it only needs to prove
$\frac{d}{dt} g(t) = 0$ along the solution curve of the matrix differential equation
\eqref{LMDE}. Actually, from equation \eqref{LMDE}, we have
\begin{align}
  g(t & +\Delta t)  = det(Q(t+\Delta t))
  = det\left((Q(t) + \dot{Q}(t) \Delta t + O\left((\Delta t)^2\right)\right)
  \nonumber \\
 &= det\left((Q(t) + SQ(t) \Delta t + O\left((\Delta t)^2\right)\right)
 = det\left(I + S \Delta t + O\left((\Delta t)^2\right)\right)det(Q(t))
 \nonumber \\
 & = \prod_{i=1}^{M}
 \left((1+\lambda_i(S) \Delta t + O\left((\Delta t)^2\right)\right)det(Q(t))
 \nonumber \\
 & = \left(1+ trace(S)\Delta t + O\left((\Delta t)^2\right)\right)g(t),
 \label{DIFGT}
\end{align}
where $\lambda_{i}(S)$ represents the eigenvalue of matrix $S$. Here, we use
the property
\begin{align}
  trace(S) = \sum_{i = 1}^{M} \lambda_{i}(S). \nonumber
\end{align}
Since matrix $S$ is skew-symmetric, we have $trace(S) = 0$. Replacing this result
into equation \eqref{DIFGT}, we have
\begin{align}
  \frac{d \, g(t)}{dt} =
  \lim_{\Delta t \to 0} \frac{g(t + \Delta t ) - g(t)}{\Delta t}
  =0, \nonumber
\end{align}
which gives the proof of the result. \eproof

\vskip 2mm

For a linear Hamiltonian dynamical system, the particle satisfies the law of
conservation of energy. Similarly, the generalized energy of the dynamical
system (\ref{LMDE}) conserves constant if we define its generalized energy
as
\begin{align}
  E(t) = trace \left(Q(t)^{T}Q(t) \right)
  = \sum_{i=1}^{M}\sum_{j=1}^{M}q_{ij}^{2}(t), \label{GENERGY}
\end{align}
where $q_{ij}(t) \ (i, \ j = 1: \ M)$ are entries of matrix $Q(t)$.

\vskip 2mm

\begin{property} \label{PROPERGEN}
The generalized energy of the dynamical system \eqref{LMDE} conserves constant.
\end{property}
\proof From the definition of the generalized energy \eqref{GENERGY},
we have
\begin{align}
  \frac{dE(t)}{dt} = \frac{d}{dt} trace\left(Q(t)^{T}Q(t)\right)
   = 2trace \left(Q(t)^{T}\dot{Q}(t) \right). \label{DGENGY}
\end{align}
Since $Q(t)$ satisfies the dynamical equation \eqref{LMDE},
replacing $\dot{Q}(t)$ with $SQ(t)$ in equation \eqref{DGENGY}, we obtain
\begin{align}
  \frac{dE(t)}{dt} = 2trace\left(Q(t)^{T}SQ(t)\right). \label{TRQSQ}
\end{align}
Noticing the trace property $trace \left(A\right) = trace \left(A^{T}\right)$
and from equation \eqref{TRQSQ}, we have
\begin{align}
  \frac{dE(t)}{dt} & = 2trace\left(Q(t)^{T}SQ(t)\right)
                    = 2trace\left(Q(t)^{T}S^{T}Q(t)\right)  \nonumber \\
                   & = -2trace\left(Q(t)^{T}SQ(t)\right) = 0, \label{RGENERGY}
\end{align}
where we use the assumption $S^{T} = -S$ in the above third equality
of equation \eqref{RGENERGY}. Thus, we prove that the generalized energy of
the dynamical system (\ref{LMDE}) conserves constant. \eproof

\vskip 2mm

Another interesting property is about the inverse proposition of property
\ref{PROPERGEN}. We state this property as the following Property
\ref{InvEngy}.

\vskip 2mm

\begin{property} \label{InvEngy}
Assume that the generalized energy $trace \left(Q(t)^{T}Q(t)\right)$ of the
dynamical system \eqref{LMDE} conserves constant, then matrix $S$ is
skew-symmetric.
\end{property}
\proof Since the generalized energy $trace\left(Q(t)^{T}Q(t)\right)$ along
the solution of equation \eqref{LMDE} conserves constant, from equations
\eqref{DGENGY}-\eqref{TRQSQ} and equation \eqref{LMDE}, we have
\begin{align}
 \frac{d \, E(t)}{dt}&= 2trace\left(Q(t)^{T}\dot{Q}(t)\right)
 = 2trace\left(Q(t)^{T}SQ(t)\right) \nonumber \\
 & = 2\sum_{i=1}^{n}\left(q(t)_{i}^{T}Sq(t)_{i}\right) = 0,
 \label{DENG}
\end{align}
where $Q(t) = \begin{bmatrix} q(t)_1, \dots, q(t)_n \end{bmatrix}$.
Let vectors $q(t)_{i} = 0 \ (i = 2, \dots, n)$ and  vector
$q(t)_{1} = e_{i} + e_{j} \, (i,\, j = 1, \dots, n)$ in equation \eqref{DENG},
we obtain
\begin{align}
  s_{ij} + s_{ji} = 0. \nonumber
\end{align}
Namely matrix $S$ is skew-symmetric. \eproof

\vskip 2mm

\begin{remark}
For the linear matrix differential equation \eqref{LMDE}, there is a
stronger property than Property \ref{PROPERGEN}. Namely $Q(t)^{T}Q(t)$ is
an invariant (see Theorem 1.6, pp. 99 in \cite{HLW2006}).
\end{remark}
\proof We denote $F(t) = Q(t)^{T}Q(t)$, then from equation \eqref{LMDE}, we
have
\begin{align}
  \frac{d \, F(t)}{dt} & = \dot{Q}(t)^{T}Q(t) + Q(t)^{T}\dot{Q}(t)
  = Q(t)^{T}S^{T}Q(t) + Q(t)^{T}SQ(t) \nonumber \\
  & = -Q(t)^{T}SQ(t)^{T} + Q(t)^{T}SQ(t) = 0, \nonumber
\end{align}
which gives the proof of the result. \eproof

\section{Pseudo-Symplectic Runge-Kutta Methods}

\vskip 2mm

Since it does not exist a general linear multiple method to satisfy the
symplectic property for a Hamiltonian dynamical system (see \cite{Tang1994}),
we consider Runge-Kutta methods with the symplecticity for linear matrix
differential equation \eqref{LMDE}. An s-stage Runge-Kutta method for equation
\eqref{LMDE} has the following form (see \cite{Butcher2003}) :
\begin{align}
   & Y_{i} = Q_{k} + \Delta t \sum_{j=1}^{s}a_{ij}SY_{j},
   \quad 1 \le i \le s, \label{MDPRK} \\
   & Q_{k+1} = Q_{k} + \Delta t \sum_{i = 1}^{s}b_{i}SY_{i}, \label{RKM}
\end{align}
where $\Delta t$ is a time-stepping length and
$b_{i},\, a_{ij} \, (i, j = 1, 2, \dots, s)$ are constants. For a Runge-Kutta
method, we can write a condensed representation
, which is so-called Butcher-array as Table \ref{RKBUCHTAB}
(see \cite{Butcher2003}).

\vskip 2mm

\begin{table}[h]
\begin{center}
 \caption{Coefficients of an RK Method.}
  \begin{tabular}{l|r}
     c & A \\
     \hline
     & $b^{T}$
  \end{tabular},
  \; $c_{i} = \sum_{j=1}^{s} a_{ij}, \, i = 1, 2, \dots, s.$
 \label{RKBUCHTAB}
\end{center}
\end{table}
The symplectic condition of a Runge-Kutta method for the even dimensional
Hamiltonian system is
\begin{align}
  M = BA + A^{T}B - bb^{T} = 0, \label{SYMCOND}
\end{align}
where $B = diag(b)$ is a diagonal matrix and $b, \, A$ are the coefficients of
the Runge-Kutta method \eqref{MDPRK}-\eqref{RKM} (see \cite{Sanz1988}, or
Theorem 4.3, p. 192 in \cite{HLW2006}, or Theorem 1.4, p. 267 in \cite{FQ2010},
or equation (6.14), p. 152 in \cite{LR2004}). For the linear matrix differential
equation \eqref{LMDE} with the odd dimensional variables, we have the same
symplectic condition \eqref{SYMCOND} of the Runge-Kutta method. We state it as
the following Theorem \ref{THELMDESYM}.
\begin{theorem} \label{THELMDESYM}
Assume that $Q_{k+1}$ is the solution of equations \eqref{MDPRK}-\eqref{RKM},
when the coefficients of a Runge-Kutta method  satisfy the symplectic condition
\eqref{SYMCOND}, then we have
  \begin{align}
     dQ_{k+1}\wedge SdQ_{k+1} = dQ_{k}\wedge SdQ_{k}. \label{DPSDSYM}
  \end{align}
\end{theorem}
\proof From equation \eqref{RKM}, we have
\begin{align}
   dQ_{k+1}\wedge SdQ_{k+1} & =
   \left(dQ_{k} + \Delta t \sum_{i = 1}^{s}b_{i}SdY_{i}\right)
   \wedge S \left(dQ_{k} + \Delta t \sum_{i = 1}^{s}b_{i}SdY_{i}\right) \nonumber \\
   & = dQ_{k}\wedge SdQ_{k} + \Delta t \sum_{i=1}^{s}b_{i} \left(dQ_{k}\wedge S^{2}dY_{i}
   +SdY_{i} \wedge SdQ_{k}\right)   \nonumber \\
   & + (\Delta t)^{2} \sum_{i=1}^{s}\sum_{j=1}^{s} \left(b_{i}b_{j} SdY_{i} \wedge S^{2}dY_{j}\right).
   \label{DQK1QK1}
\end{align}
On the other hand, from equation \eqref{MDPRK}, we have
\begin{align}
  SdY_{i}\wedge SdQ_{k} & =
  SdY_{i}\wedge S d\left(Y_{i} - \Delta t \sum_{j=1}^{s}a_{ij}SY_{j} \right)   \nonumber \\
  & = SdY_{i}\wedge SdY_{i} -
  \Delta t \sum_{j=1}^{s}a_{ij} \left(SdY_{i} \wedge S^{2}dY_{j}\right) \nonumber \\
  & = -  \Delta t \sum_{j=1}^{s}a_{ij} \left(SdY_{i} \wedge S^{2}dY_{j}\right),
  \label{DYIDQK}
\end{align}
and
\begin{align}
  dQ_{k} \wedge S^{2}dY_{i} & =
  \left(dY_{i} - \Delta t \sum_{j=1}^{s}a_{ij}SdY_{j} \right) \wedge S^{2}dY_{i} \nonumber \\
  & = S^{T}dY_{i} \wedge S dY_{i} -
  \Delta t \sum_{j=1}^{s}a_{ij}\left(SdY_{j} \wedge S^{2}dY_{i}\right) \nonumber \\
  & = -SdY_{i} \wedge SdY_{i} -
  \Delta t \sum_{j=1}^{s}a_{ij}\left(SdY_{j} \wedge S^{2}dY_{i}\right) \nonumber \\
  & = - \Delta t \sum_{j=1}^{s}a_{ij}\left(SdY_{j} \wedge S^{2}dY_{i}\right).
  \label{DQKDYI}
\end{align}
Replacing the results of equations \eqref{DYIDQK}-\eqref{DQKDYI} into equation
\eqref{DQK1QK1}, we obtain
\begin{align}
  dQ_{k+1} & \wedge SdQ_{k+1} =  dQ_{k}\wedge SdQ_{k} \nonumber \\
  & + \left(\Delta t\right)^2
  \sum_{i=1}^{s}\sum_{j=1}^{s} \left(b_{i}b_{j}- b_{i}a_{ij} - b_{j}a_{ji} \right)
  \left(SdY_{i}\wedge S^{2}dY_{j}\right). \label{SYMCONRK}
\end{align}
Thus, from equation \eqref{SYMCONRK}, we know that the result of equation
\eqref{DPSDSYM} is true if the coefficients of a Runge-Kutta method satisfy
equation \eqref{SYMCOND}. \eproof

\vskip 2mm

\begin{theorem}
Assume that the coefficients of a Runge-Kutta method satisfy the symplectic
condition \eqref{SYMCOND} and apply this Runge-Kutta method for the linear
matrix differential equation \eqref{LMDE} to obtain its numerical solution
$Q_{k}$, then we have
\begin{align}
  Q_{k+1}^{T}Q_{k+1} = Q_{k}^{T}Q_{k} = Q(t_0)^{T}Q(t_0) = I,
  \label{ORGONALITY}
\end{align}
which also gives the conservation of the discrete generalized energy
$trace\left(Q_{k}^{T}Q_{k}\right)$.
\end{theorem}
\proof From the Runge-Kutta method \eqref{RKM}, we have
\begin{align}
 Q_{k+1}^{T}Q_{k+1} = Q_{k}^{T}Q_{k} + \Delta t \sum_{i=1}^{s}
 b_{i}\left(Q_{k}^{T}SY_{i} + Y_{i}^{T}S^{T}Q_{k}\right)
 + (\Delta t)^2 \sum_{i=1}^{s}\sum_{j=1}^{s} b_{i}b_{j}Y_{i}^{T}S^{T}SY_{j}.
 \label{QK1TQK1}
\end{align}
According to equation \eqref{MDPRK}, we obtain
\begin{align}
  Q_{k}^{T}SY_{i} = Y_{i}^{T}SY_{i} - \Delta t \sum_{j=1}^{s}a_{ij}Y_{j}^{T}
  S^{T}SY_{i}, \label{QKTSYI}
\end{align}
and
\begin{align}
  Y_{i}^{T}S^{T}Q_{k} = Y_{i}^{T}S^{T}Y_{i} - \Delta t \sum_{j=1}^{s}a_{ij}Y_{i}^{T}
  S^{T}SY_{j}. \label{YITSTQK}
\end{align}
Inserting the results of equations \eqref{QKTSYI}-\eqref{QKTSYI} into equation
\eqref{QK1TQK1}, and using the symplectic condition \eqref{SYMCOND},
we have
\begin{align}
  Q_{k+1}^{T}Q_{k+1} & = Q_{k}^{T}Q_{k} + \Delta t \sum_{i=1}^{s}
 b_{i}\left(Q_{k}^{T}SY_{i} + Y_{i}^{T}S^{T}Q_{k}\right)
 + (\Delta t)^2 \sum_{i=1}^{s}\sum_{j=1}^{s} b_{i}b_{j}Y_{i}^{T}S^{T}SY_{j}
 \nonumber \\
 & = Q_{k}^{T}Q_{k} + (\Delta t)^2 \sum_{i=1}^{s} \sum_{j=1}^{s}
 (b_{i}a_{ij} + b_{j}a_{ji} - b_{i}b_{j})Y_{i}^{T}S^{2}Y_{j}
 = Q_{k}^{T}Q_{k}, \nonumber
\end{align}
which gives the proof of the result of equation \eqref{ORGONALITY} and also
gives $trace \left(Q_{k+1}^{T}Q_{k+1}\right)
 = trace \left(Q_{k}^{T}Q_{k}\right)$. \eproof

\vskip 2mm

When $s$ equals $1$ of a Runge-Kutta method (\ref{MDPRK})-(\ref{RKM}) and its
coefficients are listed by Table \ref{MIDPOMETH}, the method is also called
the implicit midpoint method with order 2. It is not difficult to verify that
the implicit midpoint method satisfies the symplectic condition \eqref{SYMCOND}.
Therefore, it is a symplectic geometric method.

\vskip 2mm

\begin{table}[h]
\begin{center}
\caption{Coefficients of the implicit midpoint Method.}
  \begin{tabular}{l|r}
     $\frac{1}{2}$ & $\frac{1}{2}$ \\
     \hline
     & 1
  \end{tabular}
  \label{MIDPOMETH}
\end{center}
\end{table}

\vskip 2mm

If we apply the implicit midpoint method to the linear matrix differential equation
\eqref{LMDE}, we have
\begin{align}
  Q_{k+1} = \left(I - \frac{1}{2}\Delta t S\right)^{-1} \left(I + \frac{1}{2}
  \Delta t S\right)Q_{k}. \label{LMDMIDR}
\end{align}
Here, the Cayley transform
\begin{align}
   \Omega_{\Delta t}(S) = \left(I-\frac{1}{2}\Delta t S\right)^{-1}
   \left(I + \frac{1}{2}\Delta t S \right) \label{CayleyTran}
\end{align}
is commutative. Namely $\Omega_{\Delta t}$  equals
$\Omega_{\Delta t}^{*}$ and the adjoint operator
$\Omega_{\Delta t}^{*}$ is defined by
\begin{align}
  \Omega_{\Delta t}^{*}(S)  =  \Omega_{-\Delta t}^{-1}(S)
   = \left(I + \frac{1}{2} \Delta t S\right)
  \left( I - \frac{1}{2} \Delta t S\right)^{-1}. \label{ACayleyTran}
\end{align}
When matrix $S$ is skew-symmetric, from equations
\eqref{CayleyTran}-\eqref{ACayleyTran}, we have
\begin{align}
  \Omega_{\Delta t}(S)^{T}\Omega_{\Delta t}(S) = I.  \nonumber
\end{align}
Therefore, from equation \eqref{LMDMIDR}, we obtain
\begin{align}
   Q_{k+1}^{T}Q_{k+1} = Q_{k}^{T}\Omega_{\Delta t}(S)^{T}\Omega_{\Delta t}(S)
   Q_{k} = I. \nonumber
\end{align}
Namely the numerical solutions of the implicit midpoint method preserve the
orthogonal invariant.

\vskip 2mm

The Cayley transform $\Omega_{\Delta t}$ is also looked as a composition of the
explicit Euler transform
\begin{align}
   \Phi_{\Delta t}(S) = (I+ \Delta t S)     \nonumber
\end{align}
and the implicit Euler transform
\begin{align}
    \Phi_{\Delta t}^{*}(S) = (I-\Delta t S)^{-1}.  \nonumber
\end{align}
That is to say
\begin{align}
   \Omega_{\Delta t}(S) =
   \Phi_{\frac{1}{2}\Delta t}^{*}(S)\Phi_{\frac{1}{2}\Delta t}(S). \nonumber
\end{align}

\vskip 2mm
\begin{definition}
The adjoint operator $\Phi_{\Delta t}^{*}$ of $\Phi_{\Delta t}$ is defined by
$\Phi_{-\Delta t}^{-1}$. If the adjoint operator $\Phi_{\Delta t}^{*}$ equals
$\Phi_{\Delta t}$, it is called symmetric.
\end{definition}
According to the definition of the symmetric operator, it is not difficult
to see that the Cayley transform \eqref{CayleyTran} is symmetric.

\section{Numerical Experiments}

\vskip 2mm

In order to illustrate the structure-preserving property of the symplectic
method for the differential equation \eqref{LMDE}, we compare the numerical
results of the symplectic implicit midpoint method listed by Table \ref{MIDPOMETH}
with the numerical results of the non-symplectic explicit second order Runge-Kutta
method listed by  Table \ref{EXP2RK} \cite{Butcher2003}.
\begin{table}[h]
\begin{center}
\caption{Coefficients of the explicit second order RK Method.}
\begin{tabular}{l|cr}
     0 &   \\
     $\frac{1}{2}$ & $\frac{1}{2}$ \\
     \hline
     & 0 & 1
  \end{tabular}
 \label{EXP2RK}
\end{center}
\end{table}
When we apply the explicit second order Runge-Kutta method to the linear
matrix differential equation \eqref{LMDE}, we obtain the following iteration
formula
\begin{align}
  Q_{k+1} = \left( I + \Delta t S + \frac{1}{2} \Delta t^2 S^2 \right)Q_{k}.
  \label{Exp2ordRK}
\end{align}
It is not difficult to verify that the coefficients of the explicit Runge-Kutta
method can not satisfy the symplectic condition \eqref{SYMCOND}. This means that
its numerical solutions of the explicit Runge-Kutta method can not preserve
geometric structure \eqref{DPSDSYM} and can not comply with the conservation of
energy of the dynamical system \eqref{LMDE}.

\vskip 2mm

The test problem is given as
\begin{align}
  S = \begin{bmatrix}
  0 & 2 & -0.1 \\
  -2 & 0 & 0 \\
  0.1 & 0 & 0
  \end{bmatrix}. \label{PROB1}
\end{align}
The integrated interval is $[0,\, 2000]$ and the fixed time-stepping length is
0.1.

\vskip 2mm

\begin{figure}[ht]
 \centerline{
  \scalebox{0.8}{\includegraphics{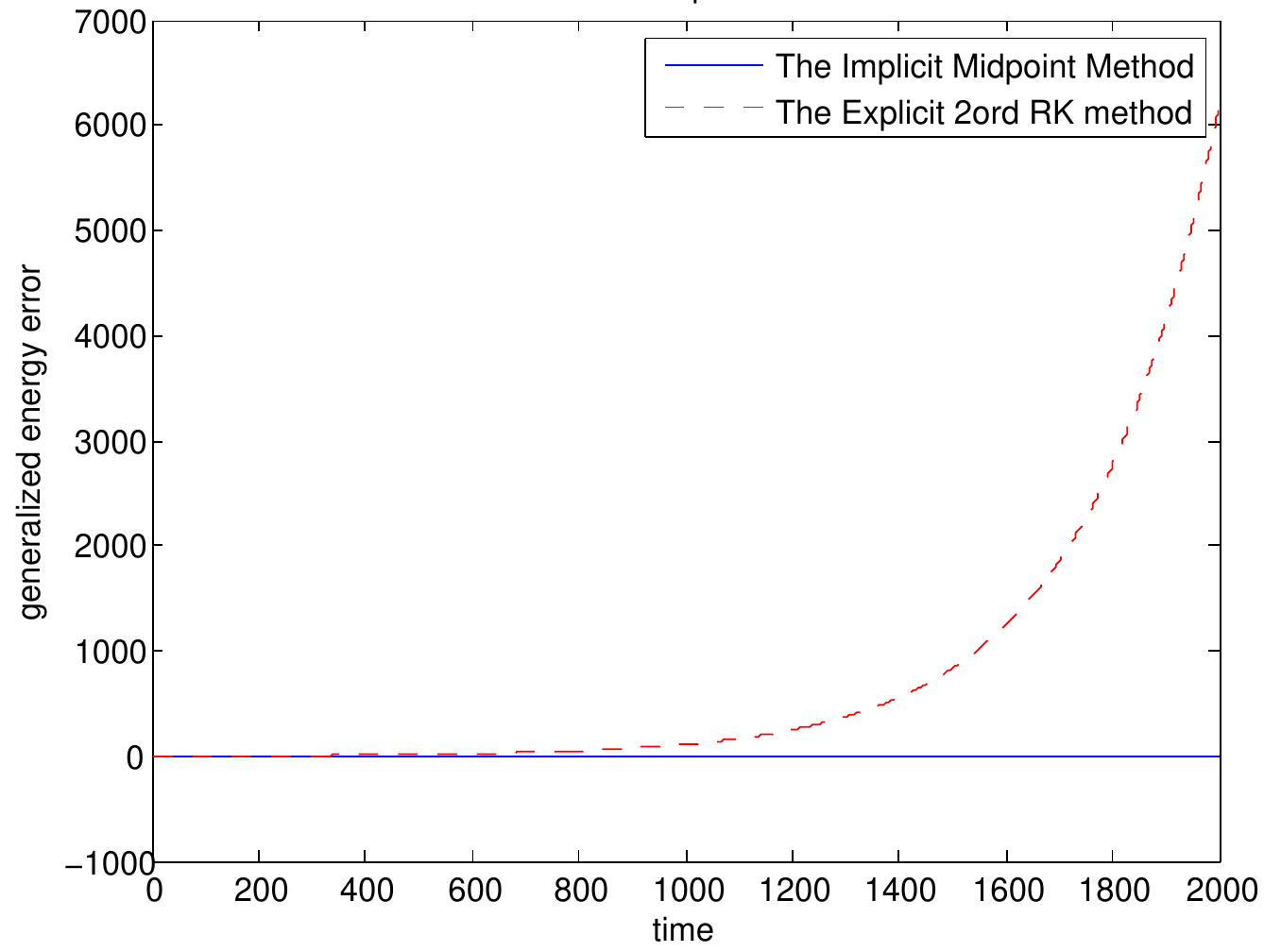}}
  }
  \caption{Numerical results of example 1} \label{Fig:Exam1}
\end{figure}

The numerical results of the test problem  are presented by Figure \ref{Fig:Exam1}.
The horizontal axis is on time and the vertical axis represents the error of the
discrete energy. From Figure \ref{Fig:Exam1}, we find that the generalized energy of
the symplectic implicit midpoint method \eqref{LMDMIDR} fluctuates tiny, and the
generalized energy of the non-symplectic explicit Runge-Kutta method grows with time. It
means that the numerical results conform to the theoretical results of the
previous sections.

\vskip 2mm

\section{Conclusion}

\vskip 2mm

We mainly extend the applicable fields of symplectic
geometric algorithms from the even dimensional Hamiltonian system to
the odd dimensional dynamical system, and discuss the geometric
and algebraic properties of symplectic Runge-Kutta methods for
the linear matrix differential equation, such as the symplecticity
and the orthogonality. It is worth noting that the implicit midpoint rule
preserves the Lie group structure of orthogonal matrices (see for
example p. 118 in \cite{HLW2006}) and this is the interested research
topic. Another interesting issue is how to preserve the invariant
of $Q_{k}^{T}Q_{k}$ when we use the approximate technique in
\cite{LW2011} if the initial matrix $Q_0$ is not orthogonal.
We would like to consider those issues in our future work.

\vskip 2mm

\section*{Financial and Ethical Disclosures}
\begin{itemize}
  \item[$\bullet$]  Funding: This study was funded by by Grant 61876199 from National
  Natural Science Foundation of China, Grant YBWL2011085 from Huawei Technologies
  Co., Ltd., and Grant YJCB2011003HI from the Innovation Research Program of Huawei
  Technologies Co., Ltd..
  \item[$\bullet$]  Conflict of Interest: The authors declare that they have no
  conflict of interest.
\end{itemize}


\end{document}